\documentclass{amsart}

\usepackage{latexsym, amssymb, amsfonts,amsmath}
\usepackage{mathtools, amsthm, amsfonts,amssymb, enumitem, mathabx}

\usepackage[
    backend=biber,
    language=english]{biblatex}
\theoremstyle{definition}

\newtheorem{Definicao}{Definition}[section]

\newtheorem*{Prova}{Proof}
\newtheorem*{Demonstracao}{Proof}
\newtheorem{Observacao}[Definicao]{Remark}

\theoremstyle{plain}

\newtheorem{Teorema}[Definicao]{Theorem}
\newtheorem{Lema}[Definicao]{Lemma}
\newtheorem{Proposicao}[Definicao]{Proposition}
\newtheorem{Corolario}[Definicao]{Corollary}

\newcommand\restr[2]{\ensuremath{\left.#1\right|_{#2}}}
\usepackage{stackengine}

\renewcommand{\L}{\ensuremath{\mathcal{L}}}

\newcommand{\F}{\ensuremath{\mathcal{F}}}

\newcommand{\M}{\ensuremath{\mathcal{M}}}

\newcommand{\Cc}{\ensuremath{\mathcal{C}}}

\newcommand{\f}{\ensuremath{\varphi}}

\renewcommand{\P}{\ensuremath{\mathcal{P}}}

\renewcommand{\S}{\ensuremath{\mathcal{S}}}

\title[Interior preorder structures and neighbourhood systems]{Translations between interior preorder structures and coherent neighbourhood systems for intuitionist modal logic}

\author{Aliel Minatti Andrade and Rogério Augusto dos Santos Fajardo}

\begin{document}

\begin{abstract}
On this article we provide a relation between two inherently different semantic structures for intuitionistic modal logic. We start by recalling the Heyting Algebras, then defining the language and the axioms for the \(iS4h\) intuitionistic calculus. We then proceed to analyse two structures discussed on \autocite{s4i} and their relation. Finally, we discuss about a structure defined on \autocite[Chapter 6]{MM}, establishing properties in order for it to be sound with \(iS4h\), and finally conclude that it is semantically equivalent to the other two structures defined by \textcite{s4i}.

\end{abstract}

\maketitle
	


\section*{Introduction}


Classical modal logic is an extension of classical propositional logic with functional symbols called \emph{modalities} and semantics based on structures of \emph{possible worlds}. The most usual approach adopts one single primitive unary modality symbol, denoted by $\Box$, and \emph{Kripke semantics}. A \emph{Kripke frame} is an ordered pair $(W, R)$, where $W$ is a non-empty set (whose elements are called \emph{possible worlds}, or simply \emph{worlds}) and $R$ is a binary relation on $W$. A formula of the form $\Box A$ is true on a possible world $w\in W$ if $A$ is true on every possible world $v$ such that $wRv$.


A system for a \emph{normal modal logic} contains, among its axioms and inference rules, propositional tautologies, substitution rule, generalization (from $A$ we deduce $\Box A$) and axiom $K$: $\Box (p\rightarrow q)\rightarrow (\Box p\rightarrow \Box q)$. 
In non-normal modal logic systems, generalization and axiom $K$ may be false.

The most general semantics for non-normal modal logic is given by \emph{minimal models} (see~\cite{chellas}, p. 207). A frame for a minimal model is a pair $(W,F)$, where $W$ is (as usual) a non-empty set and $F$ is a function from $W$ to $\mathcal{P}(\mathcal{P}(W))$. Here, a formula $\Box A$ is true on a world $w$ if, and only if, the set of the worlds where $A$ is true belongs to $F(w)$. In this way, any semantics notion for modal logic in which the truth of $\Box A$ in some world depends on the set of the worlds where $A$ is true can be translated into semantics given by minimal models. For instance, a Kripke frame $(W, R)$ can be translated to a minimal frame $(W, F)$, where $$F(w)=\left\{S\in\mathcal{P}(W) \mid\{v\in W:wRv\}\subseteq S\right\}.$$

Topology provides an interesting particular case of general frames for non-normal modal logics. Assuming we have a topology on $W$, we may define $F(w)$ as the set of all neighbourhoods of $w$. By a \emph{neighbourhood of} $w$ in a topology space $(W,\tau)$ we mean a set $S\subseteq W$ such that there is a $V\subset S$ where $w\in V$ and $V\in \tau$. I.e., a neighbourhood of $w$ is any subset of the topological space which contains an open neighbourhood of $w$. Letting $V(A)$ be the set of all possible worlds where the formula $A$ is true, straight calculations show that, in this semantics, $V(\Box A)=[V(A)]^\circ$ (the topological interior of $V(A)$).

The \emph{neighbourhood semantics} use a minimal frame $( W, F)$ which satisfies a list of axioms that is compatible to the topological interpretation given on the previous paragraph.  
The complete definitions may be found in the 6th chapter of \cite{MM} and will be presented in Definitons~\ref{vizinhanca} and~\ref{neighbourhoodmodels} of the present paper.


Until now we have mentioned only modal logics which extends classical propositional logic. I.e., we assumed so far that the propositional connectives $\wedge$, $\vee$, $\rightarrow$ and $\neg$ are interpreted as in classical propositional logic. The set $\{V(A):A\in\mathcal{L}\}$, where $\mathcal{L}$ is the modal logic language, forms a Boolean Algebra, endowed with the standard set operations. Complement is given by negation, since $V(\neg A)=W\smallsetminus V(A)$, and we have $V(A\rightarrow B)=(W\smallsetminus V(A))\cup V(B)$, which is the classical interpretation for the implication.

The focus of this paper is intuitionistic modal logic, where implication is interpreted as pseudocomplement, denoted by $\rhd$. Negation becomes a derived symbol, defining $\neg A$ as $A\rightarrow \perp$, where $\perp$ is the constant for false. Instead of a Boolean Algebra, the set $\{V(A):A\in\mathcal{L}\}$ becomes a Heyting algebra, which is widely used in Intuitionistic Logic (see~\cite{ana}).

An illustrative example of Heyting algebra is obtained by taking a topology on $W$, and defining the pseudocomplement as $A\rhd B=[(W\smallsetminus A)\cup B]^\circ$. Recalling that the complement of $A$ is defined as $A\rhd\varnothing$, we notice that, if we take $A=[0,1]\smallsetminus\{1/2\}$ in the usual topology of the line, the double complement of $A$ is $[0,1]$. Remember that, in the intuitionistic logic, $p\rightarrow\neg\neg p$ is always true but the reverse can be false. 

In order to study modal intuitionistic logic,  \cite{s4i} adopts the \emph{interior preorder semantics}, based on \emph{interior Heyting algebras}, which are a generalization of the topological example, by adding the interior operator, with properties similar to the topological ones.

The main goal of this article is to establish a relation between the interior preorder semantics  and the neighbourhood systems. We note that, in the end of \cite{s4i}, they refer to the possibility of relating their proposed semantics with neighbourhood systems.

On the first section, we dedicate ourselves towards Heyting Algebras, which are an important mathematical object for intuitionist logic semantics. On the second to fourth sections, we discuss the intuitionist calculus and semantics defined on \cite{s4i}, but emphasizing on what's most needed for the understanding of our further proposed translation. Our contributions lay on the fifth section, where we give our interpretation of the neighbourhood systems discussed on \cite{MM} and we finally establish a one-by-one correspondence between the interior preorder semantics and our version of the neighbourhood systems.



\section{Heyting Algebras}

We remind the definition of Heyting Algebras. A \emph{poset} is a pair \((X, \leq)\), where \(\leq\) is a partial order on the set \(X\). A \emph{lattice} is an algebraic structure \(L = (X, \lor, \land, \bot, \top)\) with binary associative, commutative and idempotent operations \(\lor\) and \(\land\) (namely ``joint'' and ``meet'' respectively), and constants \(\bot\) and \(\top\) which are, respectively, the neutral elements of \(\lor\) and \(\land\). Moreover, absorption rules must hold, i.e., for every \(a, b \in X\), \(a \lor (a \land b) = a\) and \(a \land (a \lor b) = a\).  

When we define the relation $\leq$ on $X$ as $a \leq b$ iff $a \land b = a$. The pair $(X,\leq)$ is a poset, which is called \textit{the induced partial order} by \(L\) on \(X\).

Also, given a set \(X\), we have that \(L = (\P(X), \cup, \cap, \varnothing, X)\) is a lattice. We refer to it as the \textit{power set lattice} of \(X\).

A \textit{Heyting Algebra} is an ordered pair \((L, \rhd)\), where \(L\) is a lattice and \(\rhd\) is a binary operation (referred as ``relative pseudo-complement'') on \(X\) such that, for every \(a, b, x \in X\):

$$ (x \land a) \leq b \iff x \leq (a \rhd b).$$

We define the \textit{pseudo-complement} of \(a\in X\) as the element \(a \rhd \bot\), which will be denoted by \(\neg a\).

\begin{Proposicao}\textup{\autocite[Proposition~2.19]{ana}}
Let \(L\) be a Heyting algebra and \(a, b, c \in X\). Then, the following holds:

\begin{itemize}
    \item[(i):] If \(a \leq b\), then \(-b \leq -a\);
    \item[(ii):] \(-\bot = \top\) and \(-\top = \bot\);
    \item[(iii):] \(a \leq -(-a)\);
    \item[(iv):] \(-(-(-a)) = -a\);
    \item[(v):] \((-a) \land (-b) = -(a \lor b);\)
    \item[(vi):] \((-a) \lor (-b) \leq -(a \land b); \) 
    \item[(vii):] \(-(a \land b) =  b \rhd (-a).\)

\end{itemize}
\end{Proposicao}

On the next section, we will define the intuitionist modal language syntactically. We shall use the knowledge of Heyting algebras for sections 3 and beyond, which are going to talk about the semantics for such language. Those semantics make heavy use of those algebras.

\section{Modal intuitionist language}

\begin{Definicao}
Let \(\L_0\) be a countable set of propositional variables. We define \(\L\) the \textit{intuitionist modal language} over \(\L_0\) with the following rules of formation:

\begin{itemize}
    \item[(1)] \(\L_0 \cup \{\bot\} \subseteq \L\);
    \item[(2)] If \(A, B \in \L\), then \(A * B \in \L\) for every \(* \in \{\lor, \land, \to\}\);
    \item[(3)] If \(A \in \L\), then \(\Box A \in \L\).
\end{itemize}
\end{Definicao}


\begin{Definicao}[Intuitionistic logic]
An intuitionistic logic (or just ``logic'', given the context of intuitionistic language) is any subset \(\Gamma \subseteq \L\) closed by Modus Ponens, uniform substitution and generalization. By generalization, we mean that \(\Box A \in \Gamma\) whenever \(A \in \Gamma\).
\end{Definicao}

We shall define here the syntactics proposed on \cite{s4i}. It is called the \textit{Hilbert \(S4\) intuitionist calculus}. Here, \(S4\) stands for transitivity and reflexivity, which are directly related to the further defined axioms \((T)\) and \((4)\). 

\begin{Definicao}[\(iS4h\) Axioms]\label{Def 6.1.4}\autocite[Definition~2.2]{s4i}
The following list consists of the \textit{iS4h axioms}. They follow for every \(p, q, r \in \L_0\).

\begin{itemize}
    \item[\((A_1)\)] \(p \to (q \to p)\);
    \item[\((A_2)\)] \(p \to (q \to r) \to ((p \to q) \to (p \to r))\);
    \item[\((A_3)\)] \(p \to (p \lor q)\);
    \item[\((A_4)\)] \(q \to (p \lor q)\);
    \item[\((A_5)\)] \((p \to r) \to ((q \to r) \to ((p \lor q) \to r)));\)
    \item[\((A_6)\)] \(p \land q \to p\);
    \item[\((A_7)\)] \(p \land q \to q\);
    \item[\((A_8)\)] \((p \to q) \to ((p \to r) \to (p \to (q \land r)));\)
    \item[\((A_9)\)] \(\bot \to p\);
    \item[\((K)\)] \(\Box(p \to q) \to (\Box p \to \Box q)\);
    \item[\((T)\)] \(\Box p \to p\);
    \item[\((4)\)] \(\Box p \to \Box \Box p\).
\end{itemize}
\end{Definicao}

\begin{Definicao}[Rules of proof]\label{Def 6.1.5}
Let \(\Gamma \subseteq \L\). The following rules of proof establish the syntactic consequences of \(\Gamma\) on \(\L\).
\begin{itemize}
\item[(Ax)] If \(A\) is an axiom of \(iS4h\), then \(\varnothing \vdash A\); 
    \item[(C)] If \(A \in \Gamma\), then \(\Gamma \vdash A\);
    \item[(MP)] If \(\Gamma \vdash A\) and \(\Gamma \vdash A \to B\), then \(\Gamma \vdash B\);
    \item[(US)] If \(\varnothing \vdash A\), then \(\varnothing \vdash A'\), where \(A'\) is any uniform substitutions of the propositional variables of \(A\);
    \item[(G)] If \(\varnothing \vdash A\), then \(\Gamma \vdash \Box A\).  
\end{itemize}

Still, if \(\varnothing \vdash A\), we say that \(A\) is a \textit{theorem} of \(\L\), and we denote \(\vdash A\). The \textit{iS4h logic} is the set that consists of every theorem of \(\L\).
\end{Definicao}

\section{Interior Heyting algebras.}
On this section we will define semantic structures for \(\L\). We emphasize that \textcite[Section 3]{s4i} is the main reference for this section.

\begin{Definicao}[Preoredered sets]
Let \(X\) be a set. A binary relation \(\leq\) on \(X\) is a  \textit{preorder} if it is reflexive and transitive. Furthermore, we say that \((X, \leq)\)\footnote{We use this as the standard notation for a preordered set. By ``standard'' we mean we will always use the same tuple unless if further noticed.} (or just \(X\)) is a \textit{preordered set} if \(\leq\) is a preorder on \(X\). 
\end{Definicao}

\begin{Definicao}
Let \((X, \leq)\) be a preordered set. We say that \(U \subseteq X\) is an \textit{upset} of \(X\) if it is closed by \(\leq\). I.e, for every \(x \in U\) and \(y \in X\) such that \(x \leq y\), we have that \(y \in U\). We denote by \(\mu(X)\) the collection of every upset of \(X\).
\end{Definicao}

\begin{Proposicao}[Heyting algebra of upsets]\textup{\autocite[Section 4.2]{bez}}
Let \((X, \leq)\) be a preordered set. Then, \(\mu(X)\) is a Heyting algebra endowed with the power set lattice operations and the relative pseudocomplement operator \(\rhd\) defined by $$a \rhd b \coloneqq \{x \in X \mid \textnormal{For every }y \in X, \textnormal{ if }x \leq y \textnormal{ and } y \in a, \textnormal{ then } y \in b  \}.$$ 
\end{Proposicao}

\begin{Observacao}\label{obs:pseudocomplement}
Considering the Heyting algebra's structure given by the previous proposition, notice that, if \(a, b \in \mu(X)\) and \(a \subseteq b\), then \(a \rhd b = X\). As a matter of fact, the implication ``if \(y \in a\), then \(y \in b\)'' is trivially satisfied for every \(y \in X\), as \(a \subseteq b\), therefore \(a \rhd b = X\).

We also notice that $a\cap (a\rhd b)\subseteq b$, for every \(a, b \in \mu(X)\). In fact, if $x\in a\cap (a\rhd b)$, it follows immediately from the definition of $a\rhd b$ and the reflexivity of $\leq$ that $x\in b$. 
\end{Observacao}

\begin{Definicao}[Interior Heyting algebra]\autocite[Definition~3.3]{s4i}\label{Def 3.5}
An \textit{interior Heyting algebra} is a Heyting algebra \(\mathcal{H} = (H, \lor, \land, \rhd, 0, 1)\) with an operator \(i \colon H \to H\) such that, for every \(a, b \in H\):

\begin{itemize}
    \item[\((I_1)\)] \(i(a) \leq a\);
    \item[\((I_2)\)] \(i(a) \leq i(i(a))\);
    \item[\((I_3)\)] \(i(a \land b) = i(a) \land i(b)\);
    \item[\((I_4)\)] \(i(1) =1\). 
\end{itemize}

\(i\) is called an \textit{interior operator} for \(\mathcal{H}\).
\end{Definicao}

\begin{Observacao}
On an interior Heyting algebra \((H,i)\), we have that \(i\) is an increasing function. Indeed, let \(a, b \in H\) such that \(a \leq b\). Then, \(a \land b = a\). So, by property \((I_3)\), it follows that \(i(a) = i(a) \land i(b)\), thus \(i(a) \leq i(b)\).
\end{Observacao}

\begin{Definicao}[Positive interior operator]\autocite[Definition~3.4]{s4i}
A \textit{positive interior operator} on a preordered set \(X\) is a function \(i \colon \mu(X) \to \mu(X)\) such that \((\mu(X), i)\) forms an interior Heyting algebra. A preordered set \((X, \leq)\) with a positive interior operator \(i\) is called an \textit{interior preorder frame} \(\F = (X, \leq, i)\)\footnote{This is the standard notation for an interior preorder frame.}
\end{Definicao}
\begin{Definicao}\label{Def 6.2.7}\autocite[Definition~3.5]{s4i}
Let \(\F \) be an interior preorder frame. An \textit{interior preorder model} for \(\F\) consists on using a function \(v \colon \L_0 \to \mu(X)\), where \(\M = (\F, v)\)\footnote{This is the standard notation for an interior preorder model.}. For every formula \(A \in \L\) we recursively define \(V \colon \L \to \mu(X)\), the extension of \(v\) on \(\L\), by:

\begin{itemize}
    \item \(\restr{V}{\L_0} = v\);
    \item \(V(\bot) = \varnothing\);
    \item \(V(A \land B) = V(A) \cap V(B)\);
    \item \(V(A \lor B) = V(A) \cup V(B)\);
    \item \(V(A \to B) = V(A) \rhd V(B)\);
     \item \(V(\Box A) = i(V(A))\).
\end{itemize}
\end{Definicao}
\begin{Definicao}\label{Def 6.2.8}
We denote \(\M, x \Vdash A\) whenever \(x \in V(A)\).
Let \(\Gamma \cup \{A\} \subseteq \L\) a set of formulas. We write \(\M, x \Vdash \Gamma\) if \(\M, x \Vdash B\) for every \(B \in \Gamma\). We say that \(A\) is an \textit{\(i\)-semantic consequence of \(\Gamma\)} if, for every interior preorder model \(\M\) and \(x \in X\), we have $$\M, x \Vdash \Gamma \implies \M, x \Vdash A.$$
We denote it by \(\Gamma \Vdash_i A\).  
\end{Definicao}

\begin{Teorema}[Soundness of \(i\)]\label{Teo 6.2.9}\textup{\autocite[Proposition~3.6]{s4i}}
Let \(\Gamma \cup \{A\} \subseteq \L\). If \(\Gamma \vdash A\), then \(\Gamma \Vdash_i A\).
\end{Teorema}

On \cite{s4i}, using rather techincal lemmas as machinery, they prove the following theorem:

\begin{Teorema}[Completeness of \(i\)]\label{Teo 6.2.10}\textup{\autocite[Theorem~4.10]{s4i}}
Let \(\Gamma \cup \{A\}\). Then,

$$\Gamma \vdash A \implies \Gamma \Vdash_i A.$$
\end{Teorema}

\section{Upset Topologies}

This section is heavily inspired on \textcite[Section 5]{s4i}. Here, we present the upset topology as an alternative to interior preorder structures, in the sense of being a structure to semantically evaluate intuitionistic modal formulas, only to later prove that those structures are semantically equivalent. We present this relation specifically to later show that, by the same approach, the neighbourhood space semantics is also equivalent to the interior preorder structures.

\begin{Definicao}\autocite[Definition 5.1]{s4i}
An \textit{upset topology} on a preordered set \((X, \leq)\) is a topology \(\tau\) on \(X\) such that \(\tau \subseteq \mu(X)\). On this case, \(\Cc = (X, \leq, \tau)\)\footnote{We use this as standard notation for up-spaces.} is said to be an \textit{up-space}. 
\end{Definicao}

Similarly to interior preorder models, it is possible to define a semantics for the language \(\L\) by using up-spaces.

\begin{Definicao}[Up-space semantics]
Let \(\Cc\) be an up-space. An \textit{up-model} \(\M = (\Cc, v)\)\footnote{This is the standard notation for an up-model} is defined with a function \(v \colon \L_0 \to \mu(X)\) called a \textit{valuation}, and its extension \(V \colon \L \to \mu(X)\) is recursively defined the same way as in \ref{Def 6.2.7}, but changing the modal clause (\(\Box\)) to: 
$$V(\Box A) = [V(A)]^\circ.$$
We denote \(\M \Vdash A\) if \(V(A) = X\), and we say that \(\M\) \textit{satisfies} \(A\).
\end{Definicao}

\begin{Definicao}
Let \(A \in \L\) and \(\M\) an up-model. For every \(x \in X\), we write \(\M, x \Vdash A\) if, and only if, \(x \in V(A)\). We say that \(A\) is \textit{valid} on \(\Cc\) if \(V(A) = X\) for every valuation on \(\Cc\), and we denote \(\Cc \Vdash A.\) Still, given a subset \(\Gamma \subseteq L\), we say that \(A\) is an \textit{up-semantic consequence of \(\Gamma\) } if \(\M \Vdash A\) for every up-model \(\M\) that satisfy every formula of \(\Gamma\).
\end{Definicao}

\begin{Lema}[Induced up-space topology]\textup{\autocite[Lemma 5.4]{s4i}}
Let \((X, \leq, i)\) be an interior preorder structure and define the following subcollection \(\tau_i \subseteq \mu(X)\):
$$\tau_i \coloneqq \{a \in \mu(X) \mid i(a) = a\}.$$ Then, \(\tau_i\) is an up-set topology for \((X, \leq)\), called the \textit{induced upset topology} by \(i\) on \(X\).
\end{Lema}

\begin{Lema}\textup{\autocite[Section 5]{s4i}}
Let \((X, \leq, \tau)\) be an up-space. Define \(i_{\tau} \colon \mu(X) \to \mu(X)\) by \(i_\tau(b) = \bigcup\{a \in \tau \mid a \subseteq b\}\), i.e, the topological interior of \(b\). Then, \(i_\tau\) is an interior operator.
\end{Lema}
\begin{Teorema}\textup{\autocite[Proposition 5.5]{s4i}}
Let \(F\) be the class of all interior preorder frames, and \(T\) the class of all up-spaces. Then, the functor:
\[\begin{array}{cccc}
   \f \colon & F & \longrightarrow & T  \\
     & (X, \leq, i) & \longmapsto & (X, \leq, \tau_i). 
\end{array}\]
Is bijective, with its inverse being:
\[\begin{array}{cccc}
   \psi \colon & T & \longrightarrow & F  \\
     & (X, \leq, \tau) & \longmapsto & (X, \leq, i_\tau). 
\end{array}\]
\end{Teorema}

\begin{Corolario}\label{Corolario 6.3.7}
Let \(\F = (X, \leq, i)\) be an interior preorder frame, \(\Cc = (X, \leq, \tau_i)\) the associated up-space and \(v \colon \L_0 \to \mu(X)\) a valuation. Then, the extensions \(V\) and \(V'\) of \(v\) on \(\F\) and \(\Cc\), respectively, coincide, i.e, \(V(A) = V'(A)\) for every \(A \in \L\).


\end{Corolario}
\begin{Prova}
Let \(A \in \L\). We prove that \(V(A) = V'(A)\) by induction on the complexity of \(A\). Once the valuations on those models coincide on propositional clauses, we shall just verify the modal step, i.e, that if \(A \in \L\) satisfies \(V(A) = V'(A)\), then \(\Box A\) also satisfies that. In fact, we have \(V'(\Box A) = [V(A)]^\circ\). On the other hand, \(V(\Box A) = i(V(A))\). But \([V(A)]^\circ = i_{\tau_i}(V(A)) = i(V(A))\), therefore \(V(\Box A) = V'(\Box A)\).
\end{Prova}

\begin{Corolario}[Soundness and Completeness]\label{Cor 6.3.8}
Let \(A \in \L\), and \(\Gamma \subseteq \L\). Then,
$$\Gamma \vdash A \iff \Gamma \Vdash_i A \iff \Gamma \Vdash_T A.$$
\end{Corolario}
\begin{Prova}
By theorems \ref{Teo 6.2.9}, \ref{Teo 6.2.10}, we have soundness and completeness of \(i\)-semantics, which yields the first equivalency symbol.  The second equivalency symbol is a direct consequence of \ref{Corolario 6.3.7}. \hfill \(\blacksquare\)
\end{Prova}

This last corollary shows that, in terms of logical theory (i.e., the formulas satisfied by each semantic), there is no difference between up-spaces and interior preorder models. This is exactly what we are going to propose next: that the neighbourhood semantics provides the exact same logical theory as those other two semantic models. 

\section{Neighbourhood Spaces}

On this section, we study \textcite[Chapter 6]{MM}, with the intention of stablishing a relation between their proposed neighbourhood spaces and the interior preorder frames studied by \textcite[Section 3]{s4i}.

\begin{Definicao}[Neighbourhood Spaces]\autocite[Section 6.4.1]{MM}\label{vizinhanca}
Let \((X, \leq)\) be a preordered set. A \textit{neighbourhood space} is a structure \(\S = (X, \leq, N)\), such that the function \[\begin{array}{cccc}
    N \colon & X & \longrightarrow & \P(\mu(X))   \\
     & x & \longmapsto & N_x
\end{array}\] satisfies \(N_x \subseteq N_y\) whenever \(x \leq y\). We will use \(\S = (X, \leq, N)\) as the standard notation for neighbourhood spaces. 


\end{Definicao}

\begin{Definicao}[Neighbourhood Models]\autocite[Section 6.4.1]{MM}\label{neighbourhoodmodels}
Let \(\S\) be a neighbourhood space. We say \(\M = (\S, v)\) is a \textit{neighbourhood model} for \(\S\), where \(v \colon \L_0 \to \mu(X)\) is a valuation. The extension \(V \colon \L \to \mu(X)\) of the model is recursively defined the same way as in \ref{Def 6.2.7}, but changing the modal clause (\(\Box\)) to the following: for every \(U \subseteq X\), we define \(i_N(U) = \{x \in X \mid U \in N_x\}\), and \(V(\Box A) = i_N(V(A))\), for every \(A \in \L\).
We write \(\M, x \Vdash A\) whenever \(x \in V(A)\).
\end{Definicao}

\begin{Observacao}
We notice that the restriction of $i_N$ to $\mu(X)$ has its range also included in $\mu(X).$  In fact, let \(U \in \mu(X)\), \(x \in i_N(U)\) and $y\in X$ such that \(x \leq y\). By Definition~\ref{vizinhanca} this implies that \(N_x \subseteq N_y\) and, hence, \(U \in N_y\), which means that $y\in i_N(U)$, proving that $i_N(U)$ is an upset. 
\end{Observacao}



\begin{Definicao}[\(\S\)-semantic consequence]
A formula \(A \in \L\) is said to be an \textit{\(\S\)-semantic consequence} of a set of formulas \(\Gamma \subseteq \L\) when, for every neighbourhood model \(\M\) for \(\S\) and \(x \in X\):
$$\M, x \Vdash \Gamma \implies \M, x \Vdash A.$$
Where the notation \(\M, x \Vdash \Gamma \) is similarly used as in Definition \ref{Def 6.2.8}. Whenever \(A\) is an \(\S\)-semantic consequence of \(\Gamma\), we write \(\Gamma \Vdash_\S A\).
\end{Definicao}


In what follows, we will prove results around which conditions \(N\) shall satisfy so that neighbourhood models satisfy axioms of the \(iS4h\) logic such as  \((K)\), \((T)\) or \((4)\). We emphasize that such specific conditions weren't shown by \textcite[Chapter 6]{MM}, as the focus of that work is not the \(iS4h\) calculus.




\begin{Definicao}[Filtered systems]
Let \(\S\) be a neighbourhood system. We say that \(\S\) is \textit{filtered} if, for every \(x \in X\) and \(U, V \in \mu(X)\): $$U, V \in N_x \iff U \cap V \in N_x.$$

Or, equivalently, if \(i_N(U) \cap i_N(V) = i_N(U \cap V)\).
\end{Definicao}


\begin{Proposicao}\label{Prop 6.4.6}
Let \(\S\) be a filtered neighbourhood system. Then,
\begin{itemize}
    \item[(i)] \((K)\) is valid on \(\S\); 

    
    \item[(ii)] \((4)\) is valid on \(\S\) if and only if \(i_N\) satisfies \((I_1)\);
    \item[(iii)] \((T)\) is valid on \(\S\) if and only if \(i_N\) satisfies \((I_2)\).
\end{itemize}
(For \((I_1)\) and \((I_2)\), recall Definition \ref{Def 3.5}).
\end{Proposicao}
\begin{Prova}

\

\begin{itemize}

    \item[(i)]

    First we notice that \(i_N\) is an increasing function. In fact, given \(a, b \in \mu(X)\) such that \(a \subseteq b\), we have \(a \cap b = a\) and, hence, \(i_N(a \cap b) = i_N(a)\). By the hypothesis of $\S$ being filtered, we have \(i_N(a) \cap i_N(b) = i_N(a)\) and, thus, \(i_N(a) \subseteq i_N(b)\) (recall that the partial order inherited by the heyting algebra \(\mu(X)\) is the set inclusion). Now, given \(\M\) an \(N\)-model for \(S\), we will prove that \(V(K)= X\). By the recursive definition of $V$ we have  \(V(K) = V(\Box(p \to q)) \rhd V(\Box p \to \Box q).\) Let \(x \in X\) and take \(y \in X\) such that \(x \leq y\) and \(y \in V(\Box(p \to q))\). We need to prove that \(y \in V(\Box p) \rhd V(\Box q)\). For that, let \(z \in X\) such that \(y \leq z\) and \(z \in V(\Box p)\). That means \(V(p) \in N_z\). Since \(y \in V(\Box(p \to q))\), we have \(V(p) \rhd V(q) = V(p \to q) \in N_y \subseteq N_z\). As \(S\) is filtered, we have \(V(p) \cap (V(p) \rhd V(q)) \in N_z\), so \(z \in i_N(V(p) \cap (V(p) \rhd V(q)))\). We showed in Remark~\ref{obs:pseudocomplement}  that \(a \cap (a \rhd b) \subseteq b\) for every \(a, b \in \mu(X)\).
     Hence, using that \(i_N\) is increasing, we conclude that \(z \in i_N(V(q)) = V(\Box q)\). Then \(y \in V(\Box p)\rhd V(\Box q) =  V(\Box p \to \Box q)\), so \(x \in V(K)\). Therefore \((K)\) is valid on \(\S\);

\item[(ii)] \((\implies) \colon\) Suppose that \((4)\) is valid on \(\S\), and let \(a \in \mu(X)\). We will prove that \(i_N\) satisfies \((I_1)\) i.e. \(i_N(a) \subseteq a\). Let \(\M\) be an \(N\)-model such that \(V(p) = a\). By definition, \(V(\Box p) = i_N(a)\). So, if \(x \in i_N(a)\), then \(\M, x \Vdash \Box p\), thus, by validity of \((4)\), as \(x \leq x\), it follows that \(\M, x \Vdash p\), then \(x \in V(p) = a\). Therefore, \(i_N(a) \subseteq a\).

\((\impliedby) \colon\) Reciprocally, suppose that \(i_N\) satisfies \((I_1)\). Let \(\M\) be an \(N\)-model. As \(V(\Box p) = i_N(V(p)) \subseteq V(p)\), we have that \(V(4) = V(\Box p) \rhd V(p) = X\), therefore \((4)\) is valid on \(\S\);

\item[(iii)] \((\implies) \colon\) Suppose that \((T)\) is valid on \(\S\). Take \(a \in \mu(X)\) and \(\M\) an \(N\)-model for \(\S\) such that \(V(p) = a\). If \(x \in i_N(a) = V(\Box p)\), then \(\M, x \Vdash \Box p\). As \(x \leq x\), by validity of \((T)\), it follows that \(\M, x \Vdash \Box \Box p\), so \(x \in V(\Box \Box p) = i_N(V(\Box p)) = i_N(i_N(a))\), then \(i_N(a) \subseteq i_N(i_N(a))\). 
Then, \(i_N\) satisfies \((I_2)\).

\((\impliedby) \colon\) Reciprocally, suppose that \(i_N\) satisfies \((I_2)\), and let \(\M\) be an \(N\)-model for \(\S\). 
As \(V(\Box p) = i_N(V(p)) \subseteq i_N(i_N(V(p)) = V(\Box \Box p)\), it follows from Remark~\ref{obs:pseudocomplement} that \(V(T) = V(\Box p) \rhd V(\Box \Box p) = X\). Therefore, \((T)\) is valid on \(\S\).
\end{itemize}
\hfill \(\blacksquare\)

\end{Prova}

Inspired by that result, we define:
\begin{Definicao}
Let \(\S\) be a filtered neighbourhood system. We say \(\S\) is \textit{coherent} if \(i_N\) satisfies the properties \((I_1)\) and \((I_2)\). 
\end{Definicao}
\begin{Observacao}\label{Obs 6.4.8}
By proposition \ref{Prop 6.4.6}, a filtered neighbourhood system \(S\) is coherent if and only if axioms \((4)\) and \((T)\) are valid on \(S\).
\end{Observacao}
\begin{Observacao}
If \(\S\) is a neighbourhood system, then \(\S\) validates axioms \(A_1-A_9\) of \(iS4h\). In fact, that occours once such axioms do not have modal connectives, and the extension \(V\) of a valuation \(v\) for \(S = (X, \leq , N)\) coincide with its extension on \(\F = (X, \leq, i)\) (where \(i\) is any interior operator) except possibly on the modal clause.
\end{Observacao}

With those tools, we may verify relations between interior preorder frames and coherent neighbourhood systems. On the following, we propose translations between such concepts.

\begin{Lema}\label{lema:coherent}
Let \(\F = (X, \leq, i)\) be an interior preorder frame. Define \(N^i \colon X \to \P(\mu(X))\) such that, for every \(x \in X\), \(N^i_x = \{a \in \mu(X) \mid x \in i(a)\}.\) Then, \(\S = (X, \leq, N)\) is a coherent neighbourhood system.
\end{Lema}
\begin{Prova}
Let us determine the operator \(i_{N^i} \colon \mu(X) \to \mu(X)\). Given \(a \in \mu(X)\), \(x \in i_{N^i}(a)\) if and only if \(a \in N^i_x\), which is equivalent to \(x \in i(a)\). Then, \(i(a) = i_{N^i}(a)\) for each \(a \in \mu(X)\), and, hence, \(i = i_{N^i}\). Since \(i\) satisfies axioms \((I_1)-(I_3)\) (by definition of interior preorder frame), it follows that \(\S\) is coherent.
\end{Prova}

\begin{Lema}
Let \(\S = (X, \leq, N)\) be a coherent neighbourhood system.  Then, \(\F = (X, \leq, i_N)\) is an interior preorder frame.
\end{Lema}
\begin{Prova}
As \(\S\) is coherent, we know that \(i_N\) satisfies properties \((I_1), (I_2)\) and \((I_3)\) then, by definition, \(\F\) is an interior preorder frame.
\end{Prova}

Finally, the following result establishes a bi-univocal correspondence between interior preorder frames and coherent neighbourhood systems.
\begin{Teorema}\label{Teo 6.4.13}
Let \(F\) be the class of interior preorder frames and \(S\) be the class of coherent neighbourhood systems. Then, the functor

\[\begin{array}{cccc}
   \sigma \colon & F & \longrightarrow & S  \\
     & (X, \leq, i) & \longmapsto & (X, \leq, N^i). 
\end{array}\]
Is bijective, and its inverse is given by:
\[\begin{array}{cccc}
   \gamma \colon & S & \longrightarrow & F  \\
     & (X, \leq, N) & \longmapsto & (X, \leq, i_N). 
\end{array}\]
\end{Teorema}
\begin{Demonstracao}
Let's prove that \(\gamma \circ \sigma = id_F\) and \(\sigma \circ \gamma = id_S\).

Let \( \F = (X, \leq, i)\) be an interior preorder frame. In the proof of Lemma~\ref{lema:coherent} we showed that $i_{N^i}=i$.  Then,
\(\gamma(\sigma(X, \leq, i)) = \gamma(X, \leq N^i) = (X, \leq, i_{N^i}) = \F\). 


On the other hand, given \(S = (X, \leq N)\) a coherent neighbourhood system, we will prove that \(N^{i_N} = N\). For that, let \(x \in X\). We know that \(a \in N^{i_N}_x\) if and only if \(x \in i_N(a)\) if and only if \(a \in N_x\), thus \(N^{i_N}_x = N_x\) for every \(x \in X\), so \(N^{i_N} = N\). Therefore, \(\sigma \circ \gamma = id_S\). \hfill \(\blacksquare\)

\end{Demonstracao}

\begin{Corolario}\label{Cor 6.4.14}
Let \(\F = (X, \leq, i)\) be an interior preorder frame and \(\S = (X, \leq, N^i)\) be the associated coherent neighbourhood system and \(v \colon \L_0 \to \mu(X)\) a valuation. Then, the extensions \(V\) and \(V'\) of \(v\) on \(\F\) and \(\S\), respectively, coincide. That is, for every \(A \in \L\), \(V(A) = V'(A)\).
\end{Corolario}
\begin{Prova}
Similarly to the proof of the corollary \ref{Corolario 6.3.7}, it's sufficient to prove that \(V(\Box A) = V'(\Box A)\) whenever \(V(A) = V'(A)\), for every \(A \in \L\). Effectively, let \(A \in \L\) be such that \(V(A) = V'(A)\). Then, \(V(\Box A) = i(V(A)) = i_{N^i}(V(A)) = V'(\Box A)\). \hfill \(\blacksquare\)
\end{Prova}

\begin{Corolario}[Soundness and completeness of coherent neighbourhood systems]
Let \(\Gamma \cup \{A\} \subseteq \L\). Then, on \(iS4h\) calculus, we have that
$$\Gamma \vdash A \iff \Gamma \Vdash_i A \iff \Gamma \Vdash_\S A.$$

\end{Corolario}
\begin{Prova}
The first equivalency symbol refers to corollary \ref{Cor 6.3.8}. The second one is a direct consequence of corollary \ref{Cor 6.4.14}. \hfill \(\blacksquare\)
\end{Prova}

\section{Final remarks}

We conclude that, by adding extra properties to the neighbourhood systems studied by \textcite[Chapter 6]{MM} (which we then called ``coherent''), semantically speaking, they come to have no difference from the interior preorder frames presented by \textcite[Section 3]{s4i}, once they both are sound and complete with respect to the \(iS4h\) calculus. As up-spaces are also semantically equivalent interior preorder frames, in terms of \(iS4h\), by transitivity, we also have that the coherent neighbourhood systems are semantically equivalent to the up-spaces.


\printbibliography








\

\author{Aliel Minatti Andrade}
\address{Instituto de Matemática, Estatística e Ciência da Computação\\
Universidade de São Paulo (USP)\\
Rua do Matão 1010, CEP 05508-090, São Paulo, SP, Brazil} 
\email{alielminatti@ime.usp.br}

\author{Rogério Augusto dos Santos Fajardo}
\address{Instituto de Matemática, Estatística e Ciência da Computação\\
Universidade de São Paulo (USP)\\
Rua do Matão 1010, CEP 05508-090, São Paulo, SP, Brazil} 
\email{fajardo@ime.usp.br}


\end{document}